\documentclass[12pt,reqno]{article}
\usepackage{amssymb}
\usepackage{amsmath}
\usepackage{amsthm}
\usepackage{amsfonts}
\usepackage{url}
\usepackage{hyperref}
\usepackage{breakurl}
\usepackage{rotating}
\usepackage[export]{adjustbox}

\def\modd#1 #2{#1\ \mbox{\rm (mod}\ #2\mbox{\rm )}}

\newcommand{\comment}[1]{}

\DeclareMathOperator{\Li}{Li}

\theoremstyle{plain}
\newtheorem{theorem}{Theorem}
\newtheorem{corollary}[theorem]{Corollary}
\newtheorem{lemma}[theorem]{Lemma}

\theoremstyle{definition}

\begin{document}
\bibliographystyle{plain}
\title{SanD primes and numbers
}
\author{Freeman J. Dyson\\
Institute of Advanced Study, \\Princeton, NJ 08540, USA\\
\href{mailto:dyson@ias.edu}{\tt dyson@ias.edu}
\and
Norman E. Frankel\\
School of Physics,\\
The University of Melbourne, Vic.\ 3010, Australia\\
\href{mailto:nef@unimelb.edu}{\tt nef@unimelb.edu.au}
\and
Anthony J. Guttmann\\
School of Mathematics and Statistics\\
The University of Melbourne\\
Vic.\ 3010, Australia\\
\href{mailto:guttmann@unimelb.edu.au}{\tt guttmann@unimelb.edu.au}
}

\date{}

\maketitle

\begin{abstract}

We define S(um)anD(ifference) numbers as ordered pairs $(p,\, q)$ such that the digital-sum $s_{10}(p\cdot q)=q-p=\Delta >0.$ We consider both the decimal and the binary cases in detail, and other bases more superficially. If both $p$ and $q$ are prime numbers, we refer to SanD {\em primes}. For SanD primes, we prove that, with one exception, notably the pair $(2,7),$ the differences $\Delta = q-p = 14+18k,\,\,k=0,1,2,\ldots.$ 

Based on probabilistic arguments, we conjecture that the number of (base-10) SanD numbers less than $x$ grows as $c_1\cdot x,$ where $c_1 = 2/3,$ while the number of (base-10) SanD primes  less than $x$ grows as $c_2\cdot x/\log^2{x},$ where $c_2 = 3/4.$ 

We calculate the number of SanD primes up to $3\cdot 10^{12},$ and use this data to investigate the convergence of estimators of the constant $c_2$ to the calculated value. Due to the quasi-fractal nature of the digital-sum function, convergence is both slow and erratic compared to the corresponding calculation for twin primes, though the numerical results are consistent with the calculated results.
\end{abstract}

\vskip .1cm

\noindent {\bf AMS Classification scheme numbers}: 11A41, 11A63, 11Y55, 11Y60

\vskip .2cm


\vskip .3cm

\section{Introduction}
\label{introduction}

In honour of the 95th birthday of one of the authors (FJD), another of the authors (NEF) coined the SanD prime problem.
S(um)anD(ifference) {\em primes} are defined to be the subset of primes $p,q \in \mathtt{PRIMES}$ with the property that $p \cdot q=r,$ where the sum of the (decimal) digits of $r,$ denoted $s_{10}(r),$ is equal to $q-p=\Delta>0.$ 

There is only one pair involving the prime 2, viz. $(2,7),$ as $2\cdot 7=14,$ and $s_{10}(14)=7-2=5.$ The next example is $(5,19),$ as  $5 \cdot 19=95,$ both 5 and 19 are primes and $s_{10}(95)=14=19-5.$
If we relax the requirement of primality, we refer to SanD {\em numbers.}

Of course the SanD numbers and SanD primes can be defined in terms of the digital sum in any base $b,$ though $b$ must be even for there to be a non-zero set of such numbers/primes (see Section \ref{arb}). Here we treat the decimal ($b=10$) and binary ($b=2$) bases in detail, and the general case more superficially. The effect of the digital sum constraint is more prominent in the decimal case. 

The study of digital sums goes back at least to Legendre  \cite{L98}. In the late 18th century he proved that 
\begin{equation}\label{L1}
s_b(n)=n-(b-1)\sum_{j \ge 1} \lfloor \frac{n}{b^j} \rfloor.
\end{equation}
 Because of the irregular nature of this function, 
 attention historically turned instead to the behaviour of the random variable $s_b(U_n),$ where $U_n$ assumes each of the values $\{0,\ldots,n-1\}$ with equal probability $1/n.$ Let $X_n=X_n(b)$ denote the random variable $s_b(U_n)$ just defined. The first asymptotic result was proved by Bush \cite{B40} in 1940, who showed that $$\mathbb{E}(X_n) \sim \frac{b-1}{2} \log_b{n}.$$ Mirsky \cite{M49} in 1949 showed that the error term in this expression is $O(1),$ a result implicit in Bush's calculation. A significant improvement was made by Delange \cite{D75} who showed that $$\mathbb{E}(X_n)- \frac{b-1}{2} \log_b{n}= F_1(\log_b{n}),$$ where $F_1(x)=F_1(x+1)$ is a continuous, periodic {\em nowhere differentiable} function. An elegant derivation of this result using the Mellin-Perron technique can be found in \cite{FGKBT}. An illuminating discussion of the properties of this function is given in \cite{CH14}, as well as an extensive bibliography and discussion of the literature on digital sums. We will not make use of this result, except in the most general sense of referring to the properties of digital sums.

One further result worthy of note is that the ordinary generating function of the digital sum $s_b(n)$ is given by Adams-Walter and Ruskejin in \cite{FAR09}, and is
$$\sum_{n \ge 0} s_b(n)z^n = \frac{1}{1-z}\sum_{m \ge 0} \frac{z^{b^m}-bz^{b^{m+1}}+(b-1)z^{(b+1)b^m}}{(1-z^{b^m})(1-z^{b^{{m+1}}})}.$$

In the next section we prove that the definition of SanD numbers and primes restricts the differences $q-p$ to a given subset of the integers. In Section 3 we study the growth in the number of SanD numbers and primes, and give probabilistic arguments that the number of decimal SanD numbers less than $x$ grows as  $(2/3)x$ as $x$ gets large, while the number of decimal SanD primes grows like $(3/4)x/\log^2{x}.$ In Section 4 we consider SanD primes with an arbitrary base $b.$ The number of such primes less than $x$ is also expected to grow as $c_b x/\log^2{x},$  and we calculate the constant $c_b.$  We show that $c_b=0$ when $b$ is odd. In Section 5 we give numerical results, notably the number of SanD primes less than $3 \cdot 10^{12},$ and show that the numerical data gives results consistent with the probabilistic arguments of the earlier section. Section 6 treats the case of binary SanD primes, which are also enumerated up to $3 \cdot 10^{12},$ and analysed. The next section gives an heuristic calculation of the number of SanD numbers less than $x$ by approximating the sum-of-digits function $s_{10}(p \cdot q)$ by an appropriately chosen Gaussian random variable. This gives rise to results in qualitative, though not quantitative agreement with the numerical data. We then compare this behaviour to that of the SanD primes.

\section{Possible values of $\Delta$ for SanD numbers and SanD primes.}
\subsection{SanD numbers}\label{sn}
{\lemma For base-10 SanD numbers, $\Delta\equiv 5$ (mod 9) or $\Delta\equiv 0$ (mod 9)}

{\proof Any natural number $n$ can be written, in decimal form, as $$n=\sum_k \alpha_k \cdot 10^k.$$ Its digital sum, $s_{10}(n)=\sum_k \alpha_k.$ Since $\alpha_k \cdot 10^k \equiv \alpha_k,$ (mod 9), working in (mod 9)  it follows that every number is equal to the sum of its digits. 

For SanD numbers we require that $s_{10}(n(n+\Delta))=\Delta.$ So $n(n+\Delta)-\Delta\equiv 0$ (mod 9) or $(n-1)(n+\Delta+1)\equiv 8$ (mod 9). This excludes the values $n+\Delta \equiv 2,\,5,\,8$ (mod 9). This leaves the values $n+\Delta \equiv 0,\,3\,,6$ (mod 9) and  $n+\Delta\equiv 1\,,4\,,7$ (mod 9). In the first case we have  $\Delta\equiv 0$ (mod 9)  and in the second case  $\Delta\equiv 5$ (mod 9). Thus possible values of $\Delta$ are $9k$ and $5+9k,$ $k=1,2,3,4,\ldots.$
{\corollary The condition $\Delta \equiv 0$ (mod 9) implies that the SanD numbers $(n,n+\Delta) \equiv (0,0)$ (mod 3).}
{\proof $n(n+\Delta)=n^2+\Delta n.$ If $\Delta \equiv 0$ (mod 9), then $\Delta \equiv 0$ (mod 3) and so $n^2 \equiv 0$ (mod 3), hence $n \equiv 0$ (mod 3).}

 {\corollary The condition $\Delta \equiv 5$ (mod 9) implies that the SanD numbers $(n,n+\Delta) \equiv (2,1)$ (mod 3).}
{\proof $(n,n+\Delta) \equiv 5$ (mod 9) so $n^2+5n \equiv 5$ (mod 9), which has solution $n \equiv 2$ (mod 3), hence $n+\Delta \equiv 1$ (mod 3).}

 \subsection{SanD primes} \label{sp}
 {\lemma For base-10 SanD  primes, $\Delta\equiv 5$ (mod 9). If $\Delta$ is odd, the only prime-pair is $(2,7).$ If $\Delta$ is even, then $\Delta=14+18k,$ with $k=0,1,2,3,4,\ldots.$}
{\proof For SanD primes we require that $s_{10}(p(p+\Delta))=\Delta.$ So $p(p+\Delta)-\Delta\equiv 0$ (mod 9) or $(p-1)(p+\Delta+1)\equiv 8$ (mod 9). This excludes the values $p+\Delta \equiv 2,\,5,\,8$ (mod 9), and since $p+\Delta$ is prime, the values $p+\Delta \equiv 3,\,6\,,9$ (mod 9) are also excluded. This leaves $p+\Delta\equiv 1\,,4\,,7$ (mod 9) giving $p \equiv 5\,,8\,,2$ respectively. In each case we have $\Delta\equiv 5$ (mod 9). If $\Delta$ is odd, the only solution is $p=2,\,p+\Delta=7,$ as for other primes $p,$ $p+\Delta$ is even. If $\Delta$ is even the only solutions are $\Delta=14+18k,$ with $k=0,1,2,3,4,\ldots.$}

{\corollary The condition $\Delta \equiv 5$ (mod 9) implies that the SanD prime pair $(p,p+\Delta) \equiv (2,1)$ (mod 3).}
{\proof For the prime pair $(2,7)$ the result is immediate by inspection. Otherwise the proof is identical to that of the preceding corollary.}

\section{ The conjectured asymptotic behaviour of SanD numbers and SanD primes}\label{snos}

In this section we give arguments, but not proofs, that the number of SanD numbers less than $x$ grows as $\frac{2}{3}x$ as $x$ gets large, while the corresponding result for SanD primes is $\frac{3x}{4\log^2{x}}.$ The absence of proofs is hardly surprising since even without the extra conditions that define SanD primes, no results for prime pairs $(p, q)$ with fixed gap $\Delta=q-p$ have been proved, despite the remarkable recent developments described in the papers of Zhang \cite{Z14} and Maynard \cite{M19}.

\subsection{SanD numbers} \label{sandnos}
Base-10 SanD {\em numbers} less than $x$ are defined as the set of ordered pairs $(a,b)$  such that $1 \le a < b \le x$ and $b-a=s_{10}(a\cdot b).$ 

There are $x(x-1)/2 \sim x^2/2$ choices for the pair $(a,b)$ such that $1 \le a < b \le x.$
The digital sum constraint implies that $s_{10}(a\cdot b)\equiv 5$ (mod 9) or $0$ (mod 9).
We conjecture that this constraint reduces the quadratic growth of number pairs to linear growth. To see this, first note that $b-a=s_{10}(a^2)$ has exactly one solution for each $a,$ namely $b=s_{10}(a^2)+a.$
So asymptotically there are precisely $x$ such numbers $\le x.$ However it is not true that $b-a=s_{10}(ab)$ has a solution $b$ for every $a,$ and it is also possible (though it occurs infrequently) that for some values of $a$ there is more than one solution $b.$ Accordingly, we write $c_1 x$ for the number of SanD numbers less than or equal to $ x$ solving $b-a=s_{10}(a\cdot b).$

A totally different, but more complicated argument is the following: In 1968 K\'atai and Mogyor\'odi \cite{KM68} proved the asymptotic normality of the sum-of-digits function with mean $M=(9/2)\log_{10}{x}$ (this was known since 1940, \cite{B40}), and variance $V=(33/4)\log_{10}(x).$ Then $s_{10}(a\cdot b) = b-a$ holds with a probability that is, for each potential pair $(a,b)$ given by the Gaussian 
\begin{equation}\label{eqnP}
P(a,b)=\frac{1}{\sqrt{2\pi V}}\exp \left ( \frac{-(b-a-M)^2}{2V} \right ).
\end{equation}
 Since both $M$ and $V$ are very small compared to $x,$ all pairs $(a,b)$ occurring with appreciable probability have $a$ and $b$ close to the square-root of $x.$ Thus $a\cdot b \sim const. x.$

From corollaries 2 and 3,  SanD numbers must satisfy $$(a,\,b)\equiv (0,0) \,\,({\rm mod}\,\,3)\,\,\,{\rm or}\,\,(2,1)\,\,({\rm mod}\,\,3).$$ Since there are nine equally likely values for $(a,\,b)$ (mod 3), this gives a probability of $2/9$ that pairs chosen at random satisfy these conditions. Choosing two numbers at random, their product is equally likely to be 0, 1 or 2 (mod 3), so each product has probability $1/3.$ The ratio of these probabilities, $2/3,$ is the constant $c_1$ above, so the number of SanD numbers is expected to grow like $2x/3.$ Numerical experimentation is consistent with this result.

\subsection {SanD primes.}\label{sandprimes}

The fact that the pair $p,\,\,p+\Delta$ are both primes suggests the (generalized) twin-prime conjecture, albeit constrained by the stringent condition on the digital sum of the product. 

As discussed, for example, by Tao in \cite{T09}, the primes are believed (not proved) to behave pseudo-randomly. This belief goes back at least to Cram\'er \cite{C36}, whose model can be easily refined, since all primes greater than 2 are odd, to one in which primes $< x$ are modelled by a set of integers such that odd integers are selected with probability $2/x.$ Further refinement of this model, as discussed in \cite{T09} leads to the  prediction that the number of twin primes $< x$ behaves as $2 C_2 \frac{x}{\log^2{x}},$ where $$C_2 = \prod_{p \ge 3\,prime}\left ( 1 - \frac{1}{(p-1)^2}\right ),$$ and is known as the Hardy-Littlewood constant \cite{HL23}. Subdominant terms  are given by the stronger conjecture that the number of twin primes $< x$ is asymptotically $2 C_2 \Li_2(x).$ Considerably greater detail is to be found in \cite{T15}.

There are known deficiencies in the refined Cram\'er model, particularly for local problems. Maier \cite{M85} obtained the (then) surprising result that the model was defective for certain short intervals between primes, while Pintz \cite{P07} showed further problems, of a global nature. Despite this,  the refined Cram\'er model does seem to predict what is believed to be the correct asymptotic behaviour of twin primes, including the Hardy-Littlewood constant \cite{HL23}.

At a similar level of assumption then, the number of unconstrained prime SanD pairs $(p, p+ \Delta) < x$  is expected to behave as $c\cdot x/\log^2{x},$ where the constant $c$ depends on $\Delta$\footnote{The dependence on $\Delta$ is irregular, depending on the prime divisors of $\Delta.$ See for example \cite{C00}. Clearly, $c(2)=C_2,$ as defined above.}. 

  In the case of SanD primes, we have shown that $\Delta = 14+18k,$ $k=0,1,2, \ldots,$ (neglecting the isolated case $\Delta=5$). However for $x=10^k,$ the number of possible choices for $\Delta$ increases roughly as $\log_{10}{k}.$ For example, for $x=10^8$ there are exactly 8 values of $\Delta$ contributing to the total number of SanD primes $< 10^8,$ as can be seen from Table \ref{tab1} below. This would imply an extra factor $\log{x}$ in the asymptotic behaviour of SanD primes. 

There is however a second constraint, which is that the digital sum must be equal to $\Delta.$ The summands of the digits of the natural numbers up to $10^n$ vary from 1 to $9\log_{10}{x},$ that is, from 1 to $9n.$ The distribution is symmetrical and unimodal. Since the number of summands is proportional to $\log{x},$ the probability of a particular summand is proportional to $1/\log{x}.$ Similarly, restricting ourselves to primes, or even twin primes, the number of summands still appears to be proportional to $\log{x},$ so the probability of a particular summand is given by the reciprocal, $1/\log{x}.$ 

Thus we see that these two effects, the infinite number of possible values for $\Delta$ and the constraint that the digital sum of the product $s_{10}(p\cdot (p+\Delta))=\Delta,$  cancel each other out. So we expect that, asymptotically, the number of SanD primes $< x$ grows as $const\cdot x/\log^2{x}.$ 

Despite the superficial similarity to twin primes discussed above, it is more appropriate to compare the SanD prime pairs with uncorrelated pairs of prime numbers.    So we will compare the number
$N_1$ of prime pairs $(a,b),$ assuming ordering $a < b,$
with $ b-a = s_{10}(a\cdot b),$ and $b < x,$ with the total number
$N_2$ of prime pairs $(a,b)$ in this range.

   We are interested in the ratio
 $$ r = N_1/N_2.$$
The number $N_2$ of uncorrelated pairs is simply
the square of the number of primes in this range.   The Prime Number
Theorem tells us that
 $$N_2 \sim (1/2) (x/\log{x})^2,$$
asymptotically for large $x$  where the factor $1/2$ comes from the ordering.

The main statistical assumption is that the ratio $r$ is a product of factors, one
for each prime divisor $q,$ with the divisibility of the candidate primes by
different divisors $q$ being uncorrelated.   For each $q,$ the factor is the
ratio of  probabilities of integer-pairs being both prime to  $q,$ with and
without the digit-sum condition.

For every prime $q$ not equal to 3, the digit-sums are distributed
randomly over all the residue classes (mod $q$).

For each of these primes, the digit-sum condition does not change the probability that an
integer-pair will both be prime to $q.$  Each of these primes
contributes a factor unity to the ratio $ r.$    Only for $q = 3$ does the digit-sum condition
change the probabilities.

Since the digit-sum is equal to $(a\cdot b)$ (mod 3), the pair $(a,b)$ must always
be  $(2,1)$
(mod 3), as proved in corollary 5.

 The chance that the elements of an uncorrelated pair $(a,b)$ are both prime to
3 is $(4/9),$ while a pair satisfying the digit-sum condition must be $(2,1)$ (mod 3) or  $(0,0)$ (mod 3), as proved in corollaries 2 and 3. Only in the first case are both prime to 3, so the probability is $(1/2).$  The factor
contributed by the prime 3 to the ratio  $ r$ is then
$$ \frac{1}{2}/\frac{4}{9} = \frac{9}{8}.$$

Multiplying all the factors together gives the result
 $$ r = \frac{9T_1}{8T_2},$$
where $T_1$  and $T_2$ are the total number of integer pairs with and without the digit-sum condition respectively. We calculated $$T_1 \sim \frac{2}{3}x,$$ the number of SanD numbers $<x,$ in subsection \ref{sandnos}, while

$$T_2 \sim \frac{x^2}{2},\,\,{\rm so} \,\,  r = \frac{3}{2x}.$$

This gives the final result, as $ x$ tends to infinity,
$$ N_1 \sim \frac{3}{2}\cdot\frac{1}{2}\frac{x}{\log^2{x}}\,\,= \frac{3x}{4\log^2{x}}.$$

\section{SanD primes with an arbitrary base.}\label{arb}
Generalising the above result to an arbitrary base, we find that for base-$b$, the number of SanD primes less than $x$ as $ x$ tends to infinity, grows as
$ c_b x/\log^2{x},$ where $$c_b=\prod_q \frac{q(q-2)}{(q-1)^2} = \prod_q \left ( 1-\frac{1}{(q-1)^2}\right ), $$ where the product is taken over prime factors $q$ of $b-1.$  (The similarity of this constant to the Hardy-Littlewood constant is noteworthy). 

This result follows from the generalisation of the statistical argument given above for the decimal case, calculating the ratio $r=N_1/N_2.$ This ratio is, as stated, a product of factors, one for each prime divisor of $q,$ with the divisibility of the candidate primes by different divisors $ q$ being uncorrelated.

It follows from Legendre's result (\ref{L1}) that the digit-sums are randomly distributed over all the residue classes (mod $q$) except for prime factors of $b-1$. (This gave $q=3$ as the only case in the decimal case $b=10$ we originally considered. Now we have the same result for base 4, as $q=3$ is the only prime factor of $b-1=3,$ while for bases 6 and 8 the only prime factors we need consider are 5 and 7 respectively. For base 16 we'd need to consider both 3 and 5).

So the probability that the elements of an uncorrelated pair $(a,b)$ are both prime to $q$ is $((q-1)/q)^2.$ We have already seen that, modulo 3, a pair satisfying the digit sum condition must be (2,1) or (0,0). Only in the first case are both prime to 3, so the relevant probability is 1/2. Now generalising this, we see that for mod 5 the relevant pairs are (0,0), (3,1), (4,2), (2,3), and only in the last three cases are both prime to 5, giving a factor 3/4. And in general this factor will be $(q-2)/(q-1).$ Thus $$r=(T_1/T_2)(q-2)/(q-1)/((q-1)/q)^2.$$ As before $T_2=x^2/2$ and $T_1=(q-1)/q,$ which follows by generalising the argument in Section \ref{snos} as follows: The probability of a randomly chosen pair satisfying the divisibility condition is $(q-1)/q^2,$ and the probability of a particular product is $1/q$, so this ratio is $(q-1)/q,$ given as 2/3 for the decimal case, where $q=3.$ Putting these factors together gives the result.
It follows that $c_b=0$ for odd bases $b.$ That is to say, there are no SanD primes in such cases. For $b=2$ one has $c_2=1.$

The analogue of Lemma 4 for base 2 SanD primes is: $\Delta=4+2k,\,\, k=0,1,2,\ldots,$ and $c_2(base_2)=1.$\\
For base 4 SanD primes it is: $\Delta=8+6k,\,\, k=0,1,2,\ldots,$ and $c_2(base_4)=3/4.$\\
For base 6 SanD primes it is: $\Delta=6+10k;\,\,8+10k, \,\, k=0,1,2,\ldots,$ and $c_2(base_6)=15/16.$\\
For base 8 SanD primes it is: $\Delta=10+14k;\,\,18+14k; \,\,20+14k,\,\, k=0,1,2,\ldots,$ and $c_2(base_8)=35/36.$\\

\section{Numerical calculation of SanD numbers and primes.}
We first wrote a Maple program to enumerate SanD primes. We wanted to provide numerical support for the conjectured behaviour, notably that the number of SanD primes $< x$ grows as $c_2x/\log^2{x}$ with $c_2=3/4.$

  In a few hours on a 4GHz Intel i7 iMac with 64Gb of memory we found all SanD primes as large as $3 \cdot 10^8,$ but convergence was irregular. Andrew Conway kindly wrote a C program that, on a larger computer with 32 cores and 256 Gb of memory enabled us to obtain SanD primes as large as $3 \cdot 10^{12}$ in a day of computing time.

In Table 1 below we give the number of SanD primes less than $x$ for various values of $x \le 3 \cdot 10^{12},$ given with the appropriate value of $\Delta.$ We have seen that with $\Delta=5$ there is only one SanD prime. With $\Delta=14$ there appears to be only 19. This is misleading. There is a large gap to the next one, which is 11000000000000003, that is, around $10^{16},$ which is beyond our enumerative ability. Indeed, for any valid value of $\Delta$ there are (probabilistically) an infinite number of SanD primes. We now sketch a constructive proof for the case $\Delta=14,$ which can be repeated {\em  mutatis mutandis} for any other valid value of $\Delta.$\\
{\bf Proof:} Assume that the primes behave like independent random variables. Consider the number $$S = 3 + 10^r + 10^s,$$ with $r,\,s > 0.$ Then
$$S(S + 14) = 51 + 2.10^{r+1} + 2.10^{s+1} + 10^{2r}+10^{2s} + 2.10^{r+s}.$$
The digital sum $(5+1+2+2+1+1+2)$ is 14 for every such product, so the number of prime-pairs is, probabilistically speaking, infinite.   QED.

Similarly, for $\Delta=32,$ for the same number $S,$ $s_{10}(S(S+32)) =32.$ For $\Delta=50,$ the appropriate choice is $S= 7 + 3\cdot 10^r + 10^s,$ with $r,\,s > 0.$ Then $s_{10}(S(S+50)) =50.$
Similar such numbers $S$ can be found for other values of $\Delta$, showing that for every valid $\Delta$ there is an infinite number of SanD numbers, and so, probabilistically speaking, an infinite number of SanD primes.

Referring again to Table 1, Richard Brent (private communication) pointed out (i) that the diagonal above which the entries are zero can be immediately predicted from the fact that $s_{10}(n) < 9d$ for $n < 10^d,$ (ii) that the maximal entry in each row occurs approximately halfway to the boundary, and (iii) that the above probabilistic argument can be extended to  conjecture the growth of $N_\Delta(x),$ the number of SanD primes $x$ with $x < X$ and difference $\Delta.$ In particular, that $N_{14}(X)\gg \log\log{X}.$

\begin{sidewaystable}[htp!]

\caption{SanD primes data. The contribution from $\Delta=5$ adds 1 to each row and is not shown here. }
\vspace{0.1cm}
\centering
\tabcolsep=0.06cm
\begin{tabular}{|c|c|c|c|c|c|c|c|c|c|c|c|c|}
\hline
$x$&$\Delta=14$&32&50&68&86&104&122&140&158&176&194\\
\hline
$10^2$&7&0&0&0&0&0&0&0&0&0&0\\
$3 \cdot 10^2$&9&4&0&0&0&0&0&0&0&0&0\\
$10^3$&11&10&0&0&0&0&0&0&0&0&0\\
$3 \cdot 10^3$&14&29&1&0&0&0&0&0&0&0&0\\
$10^4$&15&69&21&0&0&0&0&0&0&0&0\\
$3 \cdot 10^4$&16&136&109&2&0&0&0&0&0&0&0\\
$10^5$&16&218&464&14&0&0&0&0&0&0&0\\
$3 \cdot 10^5$&18&329&1310&134&0&0&0&0&0&0&0\\
$10^6$&18&451&3579&954&8&0&0&0&0&0&0\\
$3 \cdot 10^6$&19&582&7740&4099&98&0&0&0&0&0&0\\
$10^7$&19&722&15662&16417&1170&2&0&0&0&0&0\\
$3 \cdot 10^7$&19&826&27871&48714&7831&82&0&0&0&0&0\\
$ 10^8$&19&944&47206&139196&48831&1985&6&0&0&0&0\\
$3 \cdot 10^8$&19&1014&72994&315414&200810&16247&126&0&0&0&0\\
$10^9$&19&1094&106919&696450&813091&135580&3213&0&0&0&0\\
$3 \cdot 10^9$&19&1134&147652&1347257&2508310&699799&31654&88&0&0&0\\
$10^{10}$&19&1178&195617&2499225&7575349&3686127&329134&3302&0&0&0\\
$3 \cdot 10^{10}$&19&1201&247383&4213080&18918254&13982418&1995357&43223&158&0&0\\
$10^{11}$&19&1222&303418&6850021&46040607&53629221&12799997&651464&965&0&0\\
$3 \cdot 10^{11}$&19&1240&359059&10361558&97588868&163082279&56956080&5104309&18913&8&0\\
$10^{12}$&19&1247&414440&15154071&201275729&497036770&264337125&44101608&425673&911&0\\
$3 \cdot 10^{12}$&19&1262&466029&20993451&373934734&1273600647&938235422&243895420&4365872&21996&3\\
\hline
\end{tabular}
\label{tab1}
\end{sidewaystable}

Assuming that the number of SanD primes  less than $x$ grows as $ c_2\cdot x/\log^2{x}$ as argued above, we have estimated the value of the constant $c_2$ in three different ways. Firstly, as the number of primes less than $x,$ denoted as usual by $\pi(x),$ grows as $x/\log{x}$, it follows that $xT(x)/\pi(x)^2$ should converge to $c_2.$ 
This estimator is given in the third column of Table \ref{tab2}. Another estimator is $T(x)\cdot \log^2{x}/x,$ while if the asymptotics are similar to that of twin primes, $T(x)/\Li_2(x)$ would converge more rapidly. Recall that asymptotically $$\Li_2(x) = \frac{x}{\log^2{x}}\left (1+\frac{2}{\log{x}} + \frac{6}{\log^3{x}} + O\left (  \frac{1}{\log^4{x}} \right ) \right ),$$ while \cite{dlVP99}
$$ \frac{\pi^2(x)}{x}=\frac{x}{\log^2{x}}\left (1+\frac{2}{\log{x}} + \frac{5}{\log^3{x}} + O\left (  \frac{1}{\log^4{x}} \right ) \right ),$$ so these differ only in the last quoted coefficient, and even then by only 20\%.
These last two estimators are given in columns four and five of Table \ref{tab2}. Both seem to fit the SanD distribution somewhat better than the leading term, $x/\log^2{x},$ and the same is true for binary SanD primes, discussed below. This may not persist for larger values of $x$ than we are able to compute.

In no case is convergence regular, unlike the corresponding situation for primes or twin primes. This is not surprising as the SanD primes are likely to have jagged irregularities in their distribution because the digit-sum function has jagged irregularities whenever the first or second digit changes from nine to zero. 

The data in Table \ref{tab2} is totally consistent with a value of $c_2 \approx 0.75.$ Taking data for $x \ge 10^6,$ the third column entries average around $c_2=0.725,$ the fourth column average is $c_2=0.811,$ and the fifth column gives $c_2=0.721.$ This variation is indicative of the jagged convergence, and an estimate of $c_2 \approx 0.75$ seems appropriate, in agreement with our calculation above. 

\begin{table}[htp]
\caption{Decimal SanD prime analysis. $\pi(x)$ is the number of primes $< x.$ The totals include the contribution of 1 from $\Delta=5.$}
\vspace{1mm}
\centering
\tabcolsep=0.11cm
\begin{tabular}{|c|c|c|c|c|c|c|c|c|c|c|c|c|}
\hline
$x$&Total=$T(x)$&$xT(x)/\pi(x)^2$&$T(x)\log^2(x)/x$&$T(x)/\Li_2(x)$\\
\hline
$10^2$&8&1.2800&1.697&0.7804\\
$3 \cdot 10^2$&14&1.0926&1.518&0.7965\\
$10^3$&22&0.7795&1.050&0.6343\\
$3 \cdot 10^3$&45&0.7301&0.9615&0.6438\\
$10^4$&106&0.7018&0.8992&0.6533\\
$3 \cdot 10^4$&264&0.7521&0.9352&0.7161\\
$10^5$&713&0.7749&0.9450&0.7539\\
$3 \cdot 10^5$&1792&0.7954&0.9501&0.7789\\
$10^6$&5011&0.8132&0.9564&0.8021\\
$3 \cdot 10^6$&12539&0.8002&0.9297&0.7926\\
$10^7$&33993&0.7697&0.8831&0.7639\\
$3 \cdot 10^7$&85344&0.7418&0.8432&0.7375\\
$ 10^8$&238188&0.7085&0.8082&0.7141\\
$3 \cdot 10^8$&606625&0.6890&0.7704&0.6862\\
$10^9$&1756367&0.6793&0.7543&0.6770\\
$3 \cdot 10^9$&4735914&0.6809&0.7517&0.6789\\
$10^{10}$&14289952&0.6901&0.7576&0.6883\\
$3 \cdot 10^{10}$&39400953&0.6994&0.7643&0.6978\\
$10^{11}$&120276935&0.7092&0.7716&0.7078\\
$3 \cdot 10^{11}$&333472334&0.7162&0.7763&0.7149\\
$10^{12}$&1022747594&0.7231&0.7808&0.7219\\
$3 \cdot 10^{12}$&2855514856&0.7298&0.7856&0.7287\\
\hline
\end{tabular}
\label{tab2}
\end{table}

\section{Binary SanD primes.}

We have also investigated the properties of SanD primes in base 2. The number of such SanD primes $B(x)$ less than $x$ for $x=10^n, \,\,n=2,3,4,\ldots,12$ and $x=3\cdot 10^{n}$ for $n=9, \ldots, 12,$ is given in the second column of Table \ref{b2const}. Note that $B(10)=0.$

As with base-10 SanD primes, we write $B(x) \sim b_2\cdot x/\log^2{x},$ and estimate the constant $b_2$ three different ways. The results are shown in Table \ref{b2const}.
 We see that convergence is significantly smoother than in the base-10 case, but still not monotonic, due to the jagged irregularities in the digit-sum function.

Nevertheless, a glance at the table entries would suggest a limit of 1 and this is as calculated in Section \ref{arb}.
These numbers show
clearly the difference between decimal and binary digit-sums.   The
decimal sum of $x$ differs from $x$ by a multiple of 9, and this causes
the bunching of SanD primes into the groups $\Delta= 14,\, 32,\, 50,$  etc.  In
the binary case the 9 is replaced by 1, and the divisibility by 1
does not cause any bunching.  There is only the divisibility by 2
imposed by the fact that all primes after 2 are odd.   So we see
that the binary coefficients converge to the value 1 rather than
3/4.   For the binary case, there is no special prime that plays the
role of 3 in the decimal case, and every SanD
integer pair of size $x$ has an equal chance $1/\log^2{ x}$  of being a
prime-pair.

\begin{table}[htp]

\caption{Binary SanD prime analysis. $\pi(x)$ is the number of primes $< x.$}
\centering
\tabcolsep=0.11cm
\begin{tabular}{|c|c|c|c|c|c|c|c|c|c|c|c|c|}
\hline
$x$&Total=$B(x)$&$xB(x)/\pi(x)^2$&$B(x)\log^2(x)/x$&$B(x)/\Li_2(x)$\\
\hline
$10^2$&6&0.9600&1.2724& 0.5853\\
$10^3$&32&1.1338&1.5269& 0.9226\\
$10^4$&172&1.1387&1.4591& 1.0601\\
$10^5$&922&1.0021& 1.2221& 0.9749\\
$10^6$&5632&0.9140& 1.0750& 0.9016\\
$10^7$&41421&0.9378& 1.0761& 0.9308\\
$ 10^8$&335551&1.0109& 1.1386& 1.0061\\
$10^9$&2637661&1.0202& 1.1328& 1.0167\\
$3 \cdot10^{9}$&7017793&1.0090&1.1139& 1.0060\\
$10^{10}$&20619112&0.9957& 1.0932& 0.9932\\
$3 \cdot 10^{10}$&55563472&0.9863& 1.0779& 0.9840\\
$10^{11}$&167019412&0.9849& 1.0715& 0.9828\\
$3 \cdot 10^{11}$&460924135&0.9900&1.0730&0.9881\\
$10^{12}$&1410277428&0.9970& 1.0767& 0.9954\\
$3 \cdot 10^{12}$&3905976118&0.9983& 1.0747& 0.9968\\

\hline
\end{tabular}
\label{b2const}
\end{table}%

\section{Irregular convergence}
\subsection{SanD numbers}
In this section we give an heuristic calculation for the irregular behaviour of decimal SanD numbers, based on the approximation that each sum-of-digits function $s_{10}(a\cdot b)$ can be replaced by 
 a Gaussian random variable, with mean value
$M=(9/2)\log_{10}(u)$ and variance $V=(33/4)\log_{10}(u),$ where $u=a\cdot b.$  Here
(9/2) is the mean
value of a decimal digit, and (33/4) is the mean-square-deviation from the
mean, as discussed above eqn. (\ref{eqnP}).

  This approximation is good when $u$ is large and
the $\log_{10}(u)$ digits are statistically independent variables.   Then the
equation $s_{10}(a\cdot b) = b-a$ holds with a probability that is for each potential
pair $(a,b)$ equal to the Gaussian eqn. (\ref{eqnP}).

 Since $M$ and $V$ are very small
compared with $u,$ all pairs that occur with appreciable probability have $a$ and $b$ both
close to the square root of $u.$     The potential SanD numbers $(a,b)$ lie in a narrow strip around the line
$a=b.$  To accord with the SanD prime calculation, we restrict the allowed values of $b-a$ to be integers of the form $18j-4$ with $ j=1,2,3,\ldots .$
Therefore the population density of SanD numbers 
is given by the sum

    $$W(u) = \frac{1}{\sqrt{2\pi V}}†\sum_{j \ge 0} \exp\left (-162\frac{\left (j-\frac{4+M}{18}\right )^2}{V}\right ),$$
summed over integer $ j.$   The sum is strictly over positive $ j,$ but we can
extend it to all positive and negative $j$ without significant error, since
the terms with negative $j $ are much smaller than unity.

      The sum $W(u)$ can be transformed to a rapidly converging sum by using
the Poisson
Summation formula, giving
  $$W(u) = \sum_{j=-\infty}^{j=\infty} \exp\left (-\frac{V\pi^2 j^2}{162} +\pi.ij\frac{4+M}{9}\right ).$$
We keep only the three terms of the transformed sum with $  j = 0, 1$ and
$-1.$   These give
$$W(u) = 1 +2 u^{-a}\cos\left (\frac{\pi}{2}\left (\log_{10}(u)+\frac{8}{9}\right )\right ),$$
with exponent
$$a = \frac{11.\pi^2}{216.\log(10)} \approx 0.218.$$

The omitted terms with $|j|>1$ are of order $u^{-4a}$ or smaller and are
certainly negligible.
 The equation for $W(u)$  shows that the SanD numbers occur with
approximately
constant population density 1 as a function of the square-root of $u,$ with
a deviation which is a low power of $u$ multiplied by a cosine periodic in
$\log_{10}(x)$ with period 4.

        Since the digit-sums are not in fact independent random variables,
this calculation
using Gaussian probabilities is not  rigorous.

In our previous calculations, we have been counting SanD numbers and primes $(a,b)$ such that $a<b<x,$ and calculating the number of such numbers/primes $< x.$ In the above treatment, we start with the probability $P(u)$ that an integer $u$ is the product $a\cdot b$ of a SanD number pair, so typically $x=\sqrt{u}.$ 

To test this approximate treatment, we have counted SanD numbers such that $a \cdot b < u,$ for $u=10^{n/5},$ where $n=1,2,\ldots, 40.$ Denote these counts $d(n).$ For $n < 13,$ $d(n)=0.$ For $n \geq 13$ the counts $d(n)$ are
 \begin{align*}
d(n)=&1, 3, 5, 7, 10, 12, 17, 23, 27, 35, 43, 52, 62, 73, 91, 114, 141, 165,\\
& 217, 267, 334, 430, 549, 715, 902, 1143, 1442, 1782,
\end{align*}
 for $n=13,14,15,\ldots 40,$ respectively.

The probability $P(u)$ above is then predicted to be $$P(u)=\frac{1}{12\sqrt{u}}\left (1+2u^{-a}\cos \left [\frac{\pi}{2}\left (\log_{10}(u) + \frac{8}{9} \right ) \right ]\right ),$$  (The prefactor $1/(12\sqrt{u})$ is included to give the predicted asymptotic behaviour $x/6$ for the number of SanD numbers less than $x.$  The constant $1/12$ arises as we are only counting the subset of SanD numbers corresponding to $b-a=18j-4.$) 

To connect the counts $d(n)$ with this formula we require the discrete derivative of the counting function. Thus we define $$d'(n) \equiv \frac{d(n+1)-d(n-1)}{10^{(n+1)/5}-10^{(n-1)/5}}.$$

To study the fluctuations, we need to compare the calculated value based on the Gaussian approximation $$P_{fluc}(u)=P(u)-1/(12\sqrt{u})$$ with the numerical estimate obtained from the data, $$D_{fluc}(u)=d'(n)-1/(12\sqrt{10^{n/5}}).$$

We multiply both $P_{fluc}(u)$ and $D_{fluc}(u)$  by $12\sqrt{10^{n/5}},$ which makes all the fluctuations the same relative scale, and show the results in Figure \ref{fig:d20}.

\begin{figure}[h] 
   \centering
   \includegraphics[width=\textwidth]{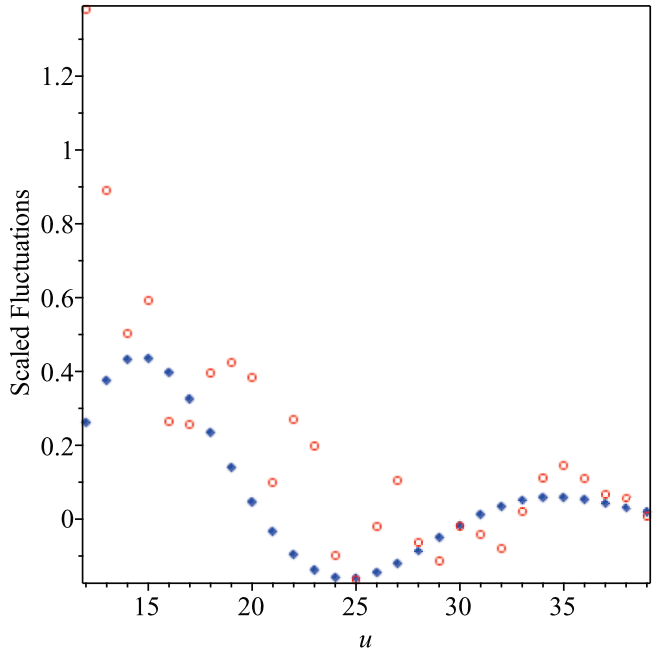} 
  \vspace{-3.3in}
 \caption{$D_{fluc}(u)\cdot 12\sqrt{10^{n/5}}$ (red circles) and $P_{fluc}(u)\cdot 12\sqrt{10^{n/5}}$ (blue diamonds) for $12 \leqslant n \leqslant 40.$  }
   \label{fig:d20}
\end{figure}

 The discrete derivatives are shown as (red) circles, the predicted probabilities as (blue) diamonds. The agreement quantitatively is disappointing, but qualitatively is instructive, showing that the actual data and the predicted data display similar irregularities. Unfortunately they don't correspond in magnitude and phase, presumably because the assumption, that the digit-sums are
independent
Gaussian random variables, is wrong.  

\subsection{SanD primes}
As discussed above, we expect the number of base-10 SanD primes $<x$ to behave as $(3/4)x/\log^2{x}$ as $x$ becomes large. To clearly see the irregular nature of convergence of the numerical data to this behaviour, we compute the deviation as follows. Recall that $T(x)$ denotes the number of SanD primes $< x.$ So $$T(x)=\frac{3}{4}\frac{x}{\log^2{x}}(1 + \theta(x)),$$ where $\theta(x)$ is of course unknown. We have calculated $\theta(x)$ for $x < 3\cdot10^{12}$ from the data in Table 2, and show the results in Figure \ref{fig:theta}. While there is insufficient data to be conclusive, there appears to be similar periodic behaviour to that observed in the SanD number fluctuations above, suggesting a possible periodic correction term.

\begin{figure}[h] 
   \centering
   \includegraphics[width=\textwidth]{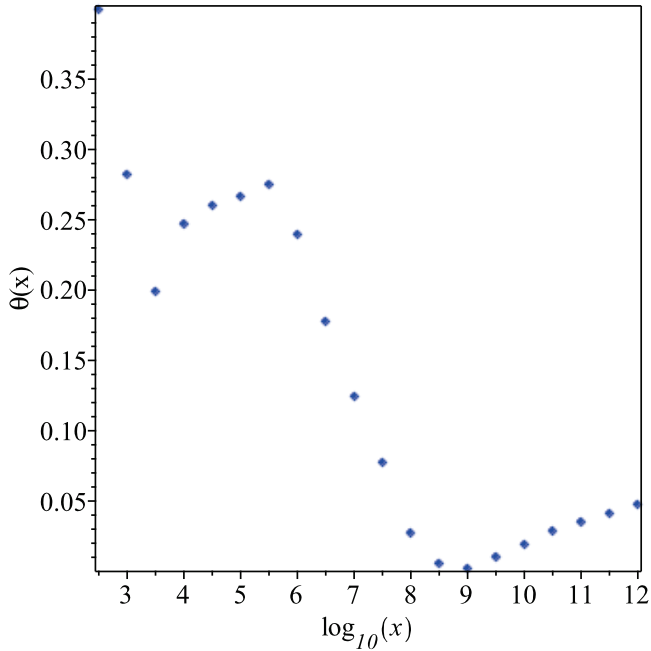} 
  \vspace{-3.4in}
 \caption{$\theta(x)$ versus $\log_{10}(x).$  }
   \label{fig:theta}
\end{figure}

\section{Conclusion}
We have defined SanD numbers as ordered pairs $(m,\, n)$ such that the digital-sum $s_{b}(m\cdot n)=n-m=\Delta > 0.$ We considered in detail both the decimal ($b=10$) and the binary ($b=2$) case. If both $m$ and $n$ are prime numbers, we refer to SanD {\em primes}. Subject to the unproven assumption that primes behave as pseudorandom numbers, in a manner described above, we show that the number of (decimal-based) SanD numbers less than $x$ grows as $c_1\cdot x,$ where $c_1 = 2/3,$ while the number of SanD primes  less than $x$ grows as $c_2\cdot x/\log^2{x},$ where $c_2 = 3/4.$ The value of the corresponding constants $c_b$ for arbitrary base-$b$ were also calculated. For binary SanD primes we show similarly that the number of such primes $B(x) < x$ behaves as $B(x) \sim b_2\cdot x/\log^2{x}$ with $b_2=1.$ 

We calculated the number of SanD numbers and primes $< 3 \cdot 10^{12}$ in order to test the above calculations. The numerical data was consistent with the conjectured results. However due to the sawtooth nature of the digital-sum function, convergence of the estimators of the constants $c_1$ and $c_2$ with increasing $x$ was found to be more erratic than the corresponding situation with twin primes, which, apart from the constant, have the same leading asymptotics. 

The twin prime distribution fits well the SanD prime pair numbers in both the decimal and binary cases (at least for primes less than $3 \cdot 10^{12}$), i.e $const\cdot \Li_2(x)$ where $const = 3/4$ and 1 respectively, in contrast with the twin primes conjecture \cite{HL23} with $const =2\cdot C_2 = 1.32....$ where $C_2$ is the twin prime constant.

\section{Acknowledgements} AJG would like to thank Andrew Conway for writing a C program to count SanD primes, Andrew Elvey Price for helpful discussions,  Richard Brent for a thorough reading of an earlier version of this paper and many suggested improvements, and Jeffrey Shallit for useful suggestions. He also gratefully acknowledges support from ACEMS the ARC Centre of Excellence for Mathematical and Statistical Frontiers.
\vspace{0.3cm}

\section{Appendix}
In the table below we show some SanD prime enumerations, giving the first 19 SanD primes for the first few values of $\Delta.$ For each entry $p$ it follows that $s_{10}(p(p+\Delta))=\Delta.$ There is one further entry, not shown, corresponding to the sole SanD prime when $\Delta=5,$ which is $p=2.$
\begin{table}[h]
\caption{Low-order SanD primes.}
\centering
\tabcolsep=0.11cm
\begin{tabular}{|c|c|c|c|c|}
\hline
$\Delta=14$&$\Delta=32$&$\Delta=50$&$\Delta=68$&$\Delta=86$\\ \hline
5 &149 &        2543 & 19961 &  412253 \\ 
17 &179 &       3137 & 28211 &  547661 \\
23& 239 &       3407 & 43541 &  871163 \\
29& 281 &       4973 & 44111 &  937661 \\
53& 389 &       5147 & 62861 &  982703 \\
59& 431 &       5693 & 66821 &  989381 \\
83& 491 &       7193 & 69941 &  992363 \\
113&509 &       7523 & 83621 &  996551 \\
167&569 &       7649 & 86561 &  999917 \\
383&659 &       7673 & 88721 &  999953 \\
443&1019 &      8243 & 89261 &  1296101 \\
1103&1031 &     8513 & 92111 &  1297601 \\
1409&1061 &     8573 & 94781 &  1329863 \\
2003&1259 &     8627 & 99191 &  1336253 \\
3203&1289 &     9293 & 120671&  1337813 \\
11483&1427 &    9461 & 125261&  1378253 \\
100043&1439 &   9497 & 129461&  1410203 \\
200003 &1901 &  9767 & 129959&  1608611 \\
1001003 & 2081 &9833& 130211&  1642211 \\

\hline
\end{tabular}
\label{tabnew}
\end{table}

\end{document}